\documentclass[12pt]{article}
\begin{document}
\def\sg{\hbox{\rm sign}\,}
\def\ds{\displaystyle}
\def\eg{{\it e.g.}\ }
\def\q{\quad} \def\qq{\qquad}
\def\l{\left} \def\r{\right}
\def\R{\hbox{\bf R}}
\def\Z{\hbox{\bf Z}}
\def\N{\hbox{\bf N}}
\def\pni{\par\noindent}
\def\vsh{\smallskip}
\def\vs{\medskip}
\def\vvs{\bigskip}
\def\vvvs{\bigskip\medskip} 
\def\vsp{\vsh\pni} 
\def\vsn{\vsh\pni}
\font\maina=cmbx12
%
  %
   %
%
   %

\centerline{{\bf FRACALMO PRE-PRINT: \ http://www.fracalmo.org}}

\hrule
\vskip 0.50truecm
\font\title=cmbx12 scaled\magstep1
\font\bfs=cmbx12 scaled\magstep1
\font\little=cmr10
\begin{center}
{\title Random walk models approximating}
 \\ [0.25truecm]
{\title symmetric space-fractional diffusion processes\footnote{Published in J. Elschner, I. Gohberg and B. Silbermann (Editors),
Problems in Mathematical Physics (Siegfried Pr\"ossdorf Memorial Volume).
Birkh¨auser Verlag, Basel/Switzerland (2001), pp. 120-145.
Series ”Operator Theory: Advances and Applications”, Vol. 121}}
\\
  [0.25truecm]
 Rudolf GORENFLO$^{(1)}$   and Francesco MAINARDI$^{(2)}$ 
\\ [0.25truecm]
$\null^{(1)}$ {\little Department of Mathematics and Informatics,
Free  University of Berlin,}
 \\ {\little Arnimallee 2-6, D-14195 Berlin, Germany}
 \\ {\little E-mail: gorenflo@mi.fu-berlin.de}
\\  [0.25truecm]
$\null^{(2)}$ {\little Department of Physics, University of Bologna, and INFN}
\\ {\little Via Irnerio 46, I-40126 Bologna, Italy}
\\{\little Corresponding Author.   E-mail: francesco.mainardi@bo.infn.it}
\end{center}

%
\vsp{\it AMS Subject Classification}:
Primary 26A33, 60E07, 60J15, 60J60;  Secondary 44A20, 45K05 
\vsp
{\it Key Words and Phrases}: 
Random walks, Riesz potential, Riesz fractional derivative, L\'evy-Feller diffusion, space-fractional diffusion, stable probability distributions, pseudo-differential operators.
%
   %
   %
%
   %
\begin{abstract}
For the symmetric case of space-fractional diffusion processes (whose
basic analytic theory has been developed in 1952 by  Feller via
inversion of Riesz potential operators) we present three random walk
models discrete in space and time. We show that for properly scaled
transition to vanishing space and time steps these models converge in
distribution to the corresponding time-parameterized stable probability
distribution. Finally, we analyze in detail a model, discrete in time but
continuous in space, recently proposed by Chechkin and Gonchar.
\\ {REMARK: } Concerning the inversion of the Riesz
potential operator $I_0^\alpha$ let us point out that its common
hyper-singular integral representation fails for $\alpha = 1$. In our
Section 2 we have shown that the corresponding hyper-singular
representation for the inverse operator $D_0^\alpha$ can be obtained also
in the critical (often excluded) case $\alpha = 1$, by analytic continuation. 
\end{abstract}


\section{Introduction: concepts and notations }
By a "space-fractional" diffusion process (or L\'evy-Feller
diffusion process) we mean a process of diffusion of an extensive
quantity with density $u(x,t)$ governed by an evolution equation
$$
{\partial u(x,t) \over \partial t} =	 D^\alpha_\theta u(x,t)
\qquad (x \in \R,\,t > 0)
\leqno(1.1) $$
with an initial condition
$$
u(x,0)=    f(x) \qquad (f\in L_1(\R)).
\leqno(1.2) $$
We interpret $x$ as space, $t$ as time variable. $D^\alpha_\theta$ is
a pseudo-differential operator acting with respect to the space
variable $x$, its symbol being
$$
\widehat{D_\theta^\alpha}(\kappa)
 = -|\kappa|^\alpha e^{\, \ds i\,(\sg \kappa)\,{\theta\pi/2 }}
 = -|\kappa|^\alpha \,i^{\,\ds (\sg \kappa)\,\theta }\,,
\leqno(1.3) $$
and the real parameters $\alpha$ and $\theta$ are restricted
by the inequalities
$$
0 < \alpha \le 2 \,, \qquad
|\theta| \le \cases{\alpha     &if $\;0 < \alpha \le 1\,,$ \cr
		    2 - \alpha &if $\;1 < \alpha \le 2\,.$\cr}
\leqno(1.4) $$
For a sufficiently well-behaved function or generalized function $\phi$
defined on $\R$ we denote by $\hat \phi$ its Fourier transform:
$$
\hat{\phi}(\kappa) = \int\limits_{-\infty}^{+\infty}
e^{i\kappa x}\, \phi(x)\,dx \qquad (\kappa \in \R)\,.
 \leqno(1.5) $$
Then,
for a generic linear pseudo-differential operator $A$
acting on these functions,
its symbol $\hat A(\kappa)$ turns out to be defined
through the Fourier representation of $\l( A\phi\r) (x)\,,$
namely
$ \widehat{A\phi}(\kappa) =
 \hat A(\kappa)\, \hat\phi(\kappa)\,.$
An often applicable practical rule is
$$
\hat A(\kappa) =    (A e^{-i\kappa x}) e^{i\kappa x}\,,
\qquad (\kappa \in \R)\,.
\leqno(1.6)
$$
If $B$ is  another pseudo-differential operator, then we have
$  \widehat{A\,B}(\kappa) =\hat A(\kappa )\, \hat B(\kappa)\,. $
Let us remark, that we chose (1.5)  to define
the Fourier transform in agreement with the common terminology of
probability theory.
\vsp
 Introducing the stable probability density
$p_\alpha(x;\theta)$ whose characteristic function (Fourier transform) is
$$
\hat p_\alpha(\kappa;\theta) =
    \exp\big(-|\kappa|^\alpha e^{\,\ds i(\sg \kappa)
{\theta\pi/2}}\big) \qquad (\kappa \in \R)
\leqno(1.7) $$
and rescaling $p_\alpha(x;\theta)$ for	$x\in \R,\ t>0$ by the similarity
variable $x\,t^{-{1\over\alpha}}$ we obtain the time-dependent stable
probability density
$$
g_\alpha(x,t;\theta) =	t^{-{1\over\alpha}}p_\alpha(xt^{-{1\over\alpha}};
\theta) \qquad (x \in \R,\, t > 0)
\leqno(1.8) $$
with which we can write the solution to
(1.1), (1.2) in the form
$$
u(x,t) = \int\limits_{-\infty}^{+\infty} g_\alpha(x - \xi,t;\theta)f(\xi)
\,d\xi.
\leqno(1.9) $$
Then, for all $t>0$, we have
$$
\cases{u(\cdot,t) \in C^\infty \cap L_1(\R)\,, \qquad
\int\limits_{-\infty}^{+\infty} u(x,t)\,dx =
\int\limits_{-\infty}^{+\infty} f(x)\,dx & \cr
f(x)\ge 0 \q \hbox{for}\q\hbox{all}\q x\in\R
\ \Rightarrow\ u(x,t)\ge 0 \q
\hbox{for}\q  \hbox{all} \q x\in\R \,. & \cr}
\leqno(1.10) $$
If the initial function is a probability density then so is also the
function $u(x,t)$ and we have in (1.1), (1.2)
the description of a Markov process.
\vsp
For orientation on the general theory of stable probability distributions
we recommend in particular [F71] and [F52], Feller's parameterization
being close to ours. In our foregoing considerations we essentially have
surveyed results of Feller's pioneering paper [F52].
For a few parameters pairs $(\alpha,\, \theta)$ representations of
$p_\alpha(x;\theta)$ in terms of elementary or
well-investigated special functions are available in the literature
but in other parameterization (see e.g. [Zo], [Sc], [SaT]).
We here content ourselves with recognizing the classical
Gauss process and the Cauchy process, respectively, in
$$
g_2(x,t;0) =	 {1\over 2\sqrt{\pi}} t^{-{1\over2}}
\exp\Big(-{x^2 \over 4t}\Big)	\,,   \qq
g_1(x,t;0) =	 {1\over \pi} {t\over x^2 + t^2}\,.
$$
\vsp
Feller in [F52] has shown that $D_\theta^\alpha$
(in case $\alpha\not = 1$)
can be viewed as inverse to the (later in [SKM] so called) Feller
potential operator which is a linear combination of two Weyl integral
operators.
In the sequel we will restrict attention to the symmetric case
$\theta=0$ retaining however, the index $0$ in order to be in concordance
with  the notation of our previous papers [GM98], [GFM], [GM99].
In [GM98]
we have discussed for $\alpha \not =1$ a random walk model for
the whole range of values $\theta$;
here we now concentrate on models not treated there, requiring however
$\theta = 0\,.$ This means that henceforth we will treat the evolution
equation
$$
{\partial u(x,t) \over \partial t} =	 D^\alpha_0 \,u(x,t)
\qquad (x \in \R,\,t > 0)
\leqno(1.11) $$
with an initial condition
$$
u(x,0)=    f(x) \qquad (f\in L_1(\R)).
\leqno(1.12) $$
The symbol of the pseudo-differential operator $D^\alpha_0$ is
$$
\widehat{D_0^\alpha}(\kappa) =	   -|\kappa|^\alpha,
\leqno(1.13) $$
and the fundamental solution to (1.11) (namely, for
$u(x,0) = \delta(x)=	$
Dirac's delta function) is the function $u(x,t) = g_\alpha
(x,t;0)$ whose Fourier transform is
$$
\widehat{g_\alpha}(\kappa,t;0) = \exp\big(-t |\kappa|^\alpha \big) \qquad
(\kappa \in \R\,, \; t>0).
\leqno(1.14) $$
Then, see (1.7) and (1.8), we get  by the Fourier inversion formula
$$
g_\alpha(x,t;0) = {1\over 2\pi}\int_{-\infty}^{+\infty}
  e^{-i\kappa x}\,\exp( -t|\kappa|^\alpha) \,d\kappa =
t^{-{1\over \alpha}}\,p_\alpha(xt^{-{1\over \alpha}};0)\,.
\leqno(1.15) $$
We will denote by $S(t)$ for $t>0$ the random variable whose
probability density is given by   $g_\alpha(x,t;0)\,.$
\vsp
{\bf Remark 1.1.}
We call a process described by (1.1), (1.2) a
"space fractional" diffusion process because it is a "fractional"
generalization of the classical diffusion process which is recovered by
taking $\alpha =2\,,$ $\theta =0\,. $
Writing (1.13) in the form
 $\widehat{D_0^\alpha}(\kappa)	= - (\kappa^2)^{\alpha /2}$
and    observing that the operator $D^2$, defined by
$\l(D^2 \,\phi\r) (x) = {d^2 \phi(x)\over dx^2}\,,$
has the symbol $-\kappa ^2\,, $ we see that
$D_0^ \alpha  = - \l(-D^2\r)^{\alpha /2}\,, $
hence the operator  $D_0^ \alpha$ is the negative of a fractional
power of the (positive definite) operator $-D^2\,. $
 By  calling a process described by (1.1), (1.2) also a
L\'evy-Feller diffusion process  we honour both L\'evy
and Feller for their essential contributions
[L25], [L54], [F52].

\vskip 0.5truecm \noindent
Our aim in the following sections is to derive
discrete-space discrete-time
random walk models approximating the
space-fractional diffusion process (henceforth considered as
a Markov process) described by (1.11) with
$u(x,0)=\delta(x)\,.$
We shall show that for properly scaled transition to
vanishing space and time steps there is convergence in distribution to
the probability distribution whose density is $g_\alpha(x,t;0)$.
We shall  give heuristic motivations for choosing our concrete models,
using in a formal way calculations with symbols.
The lack of rigour in these derivations
will hopefully not be too annoying to the pure analyst,
it will be remedied in the final end by rigorous proofs of convergence.


\section{Operators and symbols}
In this section we give a survey on the relevant operators and their
symbols thereby always assuming
$$
0<\alpha\le 2.
\leqno(2.1) $$
For general orientation and more rigorous treatment we refer to [SKM],
[R] and [F52]. We need the operator $D$ of differentiation, the Weyl
operators $I_+^\alpha$, $I_-^\alpha$ and their (formal) inverses
$I_+^{-\alpha}$, $I_-^{-\alpha}$, the Riesz potential operator
$I_0^\alpha$ whose negative inverse $-I_0^{-\alpha}$ is our
pseudo-differential operator $D_0^\alpha$ (if $\alpha\not =	1$), and
the  Hilbert
transform operator $H$. For sufficiently well behaved functions $\phi$
defined on $\R$ and if required with appropriate understanding of the
occurring integrals as Cauchy principal values we have (for $x\in \R$)
$$
(D\,\phi)(x) =	   {d\over dx}\phi(x) =    \phi'(x),
\leqno(2.2) $$
$$
\cases{
(I_+^\alpha\,\phi)(x) =    {\ds {1 \over \Gamma(\alpha)}}
 \int_{-\infty}^x (x-\xi)^{\alpha-1}\,\phi(\xi)\,d\xi\,,&\cr \cr
(I_-^\alpha\,\phi)(x) = {\ds {1 \over \Gamma(\alpha)}}
 \int_x^{+\infty} (\xi-x)^{\alpha-1}\,\phi(\xi)\,d\xi\,,& \cr}
\leqno(2.3) $$
$$
I_\pm^{-\alpha} =     \cases{
 \pm D I_\pm^{1-\alpha} &if $\; 0 < \alpha < 1\,,$    \cr\cr
{D^2}I_\pm^{2-\alpha} &if $\; 1 < \alpha \le 2\,,$ \cr}
\leqno(2.4) $$
$$
(I_0^{\alpha}\,\phi)(x)=
    {\ds{(I_+^\alpha \,\phi)(x)+(I_-^\alpha\,\phi)(x)
     \over 2\cos {(\alpha \pi/ 2)}}} =
 {\ds{\int_{-\infty}^{+\infty}|x -\xi|^{\alpha -1}\, \phi(\xi)\, d\xi
\over 2\Gamma(\alpha)\,\cos{(\alpha\pi/2)}}} \q
\hbox{for} \q\alpha  \not =	1,
\leqno(2.5) $$
$$
(H\phi)(x)={1\over\pi} \int_{-\infty}^{+\infty}
{\phi(\xi)\over x-\xi}\, d\xi\,.
 \leqno(2.6) $$
\vsp
Recalling formula (1.6) for calculation of symbols we get
by direct computation
$$
\hat D=  -i\kappa\,,
\leqno(2.7) $$
$$
\hat H =     i \,\hbox{sign}\, \kappa\,.
\leqno(2.8) $$
We take the symbols of the Weyl operators and the Riesz potential
operators from [R, Theorem 4.10] as
$$
\widehat{I_{\pm}^\alpha} (\kappa) =	(\mp i\kappa)^{-\alpha} =
|\kappa|^{-\alpha} e^{\,\ds\pm i\,(\sg \kappa)\,{\alpha\pi/ 2}}=
|\kappa|^{-\alpha}\, i^{\,\ds \pm\alpha\,\sg \kappa}\,,
\leqno(2.9) $$
from which by addition we get
$$
\widehat{I_0^\alpha} (\kappa) =    |\kappa|^{-\alpha}\,.
\leqno(2.10) $$
As already remarked and as used in [F52]
$$
{D_0^\alpha}=-I_0^{-\alpha}\q \mbox{for}\q \alpha\not = 1
\leqno(2.11) $$
in agreement with the property of symbols
$$
\widehat{D_0^\alpha}(\kappa) =	   -|\kappa|^\alpha =
-\left(\widehat{I_0^\alpha}(\kappa)\right)^{-1}\q \mbox{for}
\q \alpha\not = 1\,.
\leqno(2.12) $$
In the special case $\alpha=1$ we observe that by (2.7) and (2.8)
$$
\widehat{D_0^1}(\kappa)= -|\kappa|= -\hat D(\kappa)\, \hat H (\kappa)\,,
\leqno(2.13) $$
hence as already observed in [F52] (with a sign-modified version of the
Hilbert transform)
$$
D_0^1= - D\, H\,.
\leqno(2.14) $$
Let us yet exhibit another representation of our pseudo-differential
operator $D_0^\alpha\,.$ Via the semi-group property (see [F52] or [SKM])
$$
I_0^\alpha \,I_0^\beta = I_0^{\alpha+\beta} \qq \hbox{if}
\q \alpha,\beta \in (0,1)\,,\q \alpha + \beta < 1\,,
$$
analytic continuation to negative exponents can be justified and thus
from (2.5) and (2.11) the formula
$$
D_0^\alpha =  -{1\over 2\cos {(\alpha\pi/2)}}(I_+^{-\alpha}
+I_-^{-\alpha})\qq \mbox{for}\q \alpha \not = 1\,.
\leqno(2.15) $$
In Section 4 we shall find random walk schemes in the case
$\alpha \not =1$
by approximating in (2.15) the operators $I_+^{-\alpha}$ and
$I_-^{-\alpha}$ with the Gr\"unwald-Letnikov discretization.
In Section 5 a second
random walk scheme, for the whole range $0<\alpha\le 2$, will be obtained
by a straightforward discrete approximation of hypersingular integrals
for $I_0^{-\alpha}$ and $DH$. From [SKM, formula (12.1')] we take
(for $0<\alpha<2\,,\ \alpha\not =  1$)
$$
(I_0^{-\alpha} \phi )(x) =
     {1\over 2\Gamma(-\alpha)\,\cos {(\alpha \pi/ 2)}}
\int_0^\infty  \!\!
{\phi(x+\xi)-2\phi(x)+\phi(x-\xi)\over \xi^{\alpha +1}}\,d\xi.
\leqno(2.16) $$
\vsp
Quite formally we can obtain (2.16) by replacing in (2.5)
the integrand  by $|\xi|^{\alpha -1}\, \phi(x-\xi)\,,$
then replacing $\alpha$ by $-\alpha$, splitting
$\int_{-\infty}^\infty =
\int_{-\infty}^0 + \int_{0}^\infty$, here regularizing the right hand
side hypersingular integrals by subtracting $\phi (x)$ in the numerators,
finally substituting $-\xi$ for $\xi$ in the first right hand side
integral and then putting both integrals together.
\vsp
For convenience we simplify the coefficient in (2.16), introducing
$$
b(\alpha): =- {1\over 2\Gamma(-\alpha)\,\cos {(\alpha \pi/ 2)}}
= {1\over \pi}\, \Gamma(\alpha+1)\, \sin {(\alpha \pi/2)}\,,
\leqno(2.17)
$$
where for the latter equality we have used
$\,\sin(\alpha \pi)= 2\,\sin{(\alpha\pi/ 2)}\,\cos{(\alpha\pi/2)}\, $
and the reflection formula for the gamma function
$\Gamma(-\alpha)\, \Gamma(\alpha+1) = -\pi/\sin(\alpha\pi)\,.$
Then, for $ 0<\alpha<2,\ \alpha\not = 1\,,$
$$
(D_0^\alpha \phi)(x)=	 -(I_0^{-\alpha} \phi )(x) =
   b(\alpha)\int_0^\infty
{\phi(x+\xi)-2\phi(x)+\phi(x-\xi)\over \xi^{\alpha+1}}\,d\xi  \,.
\leqno(2.18) $$
\vsp
Note that for $\phi\in C^2(\R)$ and $\phi(x)$ bounded the integral is
finite as an improper Riemann integral and observe that $b(\alpha)>0$ for
the admitted values of $\alpha$.
\vsp
Because (2.17) gives $b(1)= 1/\pi$ we are tempted to believe the formula
$$
(D_0^1 \phi)(x)=    {1\over \pi}\int_0^\infty
{\phi(x+\xi)-2\phi(x)+\phi(x-\xi)\over \xi^2}\,d\xi\,.
\leqno(2.19) $$
We will use (2.18) and (2.19)
as  motivation for our second random walk
scheme in which  the parameter value $\alpha = 1$
does no longer play a special role. We can indeed obtain (2.19)
by looking at (2.14), then formally differentiating (behind the
integral sign) (2.6) and then splitting and regularizing the
resulting hypersingular integral in the same way as we have done in the
case $\alpha \not = 1\,.$
\vsp
In view of (2.18) and (2.19) we now have,
with $b(\alpha)$ given in (2.17),
$$
(D_0^\alpha \phi)(x)=	  b(\alpha)\int_0^\infty
{\phi(x+\xi)-2\phi(x)+\phi(x-\xi)\over \xi^{\alpha +1}}\,d\xi \q
\mbox{for}\q 0<\alpha<2\,.
\leqno(2.20) $$
Unfortunately, this formula loses its meaning in the case $\alpha =2\,.$


\section{General structure of the random walk models}
What do we mean by a {\it random walk model,  discrete in space
and   discrete in time, for a Markov process}?
Let us be given a random variable $Y$ taking
its values in the set $\Z$ of integers according to the probabilities
$$
P(Y=	k) =
p_k \qquad (k \in \Z)
\leqno(3.1) $$
with
$$
\hbox{all}\ \ p_k\ge 0\ \ \hbox{and}\ \ \sum_{k\in \Z} p_k =1\,.
\leqno(3.2) $$
We discretize the space variable $x\in \R$ and the time variable $t\ge 0$
by grid points $x_j= jh$ and instants
$t_n=  n\tau$, with $h>0,\ \tau>0, \ j\in \Z,\ n\in \N_0$.
Then, defining random variables
$$
S_n= \sum_{m=1}^n \l(h\,Y_m \r)= h\,\sum_{m= 1}^n Y_m \qq (n\in \N)
 \leqno(3.3) $$
with the $Y_m$ as {\it independent identically distributed}
 random variables,
all having the same probability distribution as the random variable $Y$,
we interpret $S_n$ as the position at time $t_n$ of a random walker
starting in the point $x= 0$ at $t= 0\,.$
Denoting by $y_j(t_n)$ the probability of sojourn of the walker in
point $x_j$ at instant $t_n$, the recursion $S_{n+1} = S_n + hY_n$
implies
$$
y_j(t_{n+1}) =	\sum_{k\in{\Z}}
p_k\,y_{j-k}(t_n) \qquad (j \in {\Z},\, n \in {\N}_0),
\leqno(3.4) $$
and the walker starting at point $x_0 = 0$ means $y_0(0) = 1$ and
$y_j(0) =  0$ for $j \ne 0\,.$ However, in the recursion scheme
(3.4) it is legitimate to use a more general initial sojourn probability
distribution ($y_j(0)| \, j \in {\Z}$).
\vsp
There is yet another possible interpretation of (3.4), namely as a
{\it scheme of redistribution of an extensive quantity}
 (e.g. mass, charge, in the random walk interpretation probability),
$y_j(t_n)$ being considered
as a clump of this extensive quantity sitting in point $x_j$ at instant
$t_n\,.$ Then (3.4) describes a  conservative and non-negativity
preserving redistribution scheme. In fact, for all $n \in {\N}$ it
follows from (3.2) that, in analogy to (1.10),
$$
\sum_{j\in {\Z}} y_j(t_n) =	\sum_{j\in {\Z}}
y_j(0) \qquad\hbox{if }\ \ \sum_{j\in {\Z}} |y_j(0)| < \infty,
$$
$$
\hbox{all }\  y_j(t_n) \ge 0 \qquad \hbox{if all} \q y_j(0) \ge 0.
$$
Such discrete redistribution schemes have been used by one
of the authors in discretization of  diffusion processes governed by
second order linear parabolic differential equations ([G70], [G78],
[GN]) as they discretely imitate essential properties of the continuous
process.
\vsp
We come nearer to the Cauchy problem (1.11), (1.12)
by intending
$y_j(t_n)$ as approximation to
$$
\int\limits_{x_j-{h/ 2}}^{x_j+{h/ 2}} u(x,t_n)\,dx
$$
which, if $u(\cdot,t_n)$ is continuous, is $\approx hu(x_j,t_n)$ for
small $h$.
\vsp
We want to show that for proper choice of the probability distribution
of the random variable $Y$ and well-scaled transition
$$
\tau =	\sigma(h)\,,\q \sigma \q \hbox{strictly monotonic},
   \q \sigma(h)\to 0 \q \hbox{as}\q h\to 0
\leqno(3.5) $$
the random walk "converges" in some sense to the Markov process
with density
$u(x,t)$ described by (1.11), (1.12) in the case that the
initial function is a probability density. More specifically, we will
prove for {\it fixed} $t>0$, $t= t_n= n\tau$ with
$\N \ni n=  t/\tau \to\infty$
(and proper scaling of $h$ and $\tau$)
that the random variable
$S_n$ of (3.3) converges {\it in distribution}
(other terminology: {\it in law}) to the random variable
$S(t)$ whose density is $g_\alpha(\cdot,t;0)\,,$
 the fundamental solution  (1.15) of (1.11).
 Observing that the
distribution function
$\,G_\alpha(x,t;0) = \int_{-\infty}^x g_\alpha(\xi,t;0)\, d\xi \,$
is continuous in $x$ (due to the fast decay in $|\kappa|$ of
$\hat g_\alpha(\kappa,t;0) =  \exp(-t|\kappa|^\alpha)$ the density
$g_\alpha(\cdot,t;0)$ is in $C^\infty(\R)$) and invoking the continuity
theorem of probability theory (see, e.g., [B, Theorem 8.28]), all we have
to do is to show that for all $\kappa \in \R$
the characteristic function
$\hat y(\kappa,t;h)$ of the random variable $S_n$ tends to
$\exp(-t|\kappa|^\alpha)$ as $h\to 0\,.$
 Note the equivalences following
from  $t=    t_n =    n\tau$ and (3.5)
$$
n\to \infty \,\Leftrightarrow\, h\to 0 \,\Leftrightarrow\,\tau \to 0\,.
\leqno(3.6) $$
\vsp
The general form of the characteristic function $\hat y(\kappa,t;h)$
can be found via the generating functions
$$
\tilde p(z) =  \sum_{j\in {\Z}}\,p_j z^j, \qquad
\tilde y(z,t_n) = \sum_{j\in {\Z}} y_j(t_n)\,z^j.
\leqno(3.7) $$
As probabilities both the $p_j$ and $y_j(t_n)$ sum up to 1 if added over
the index $j$, hence these series converge absolutely and uniformly on
the periphery $|z| = 1$ of the unit circle, and so
the functions $\tilde p$
and $\tilde y_n$ are there uniformly continuous.
The random walk $S_n$
starting at $x= 0$, we have (using
the Kronecker symbol) $y_j(0)= \delta_{j0}$, and the recursion
(3.4) being a discrete convolution we get
$$
\tilde y(z) =	 (\tilde p(z))^n.
\leqno(3.8) $$
Replacing in (3.7) $z^j$ by $e^{i\kappa x_j}= e^{i\kappa jh}$ we
obtain the corresponding characteristic functions ($\kappa \in \R$)
$$
\hat p(\kappa;h)=    \tilde p(e^{\,\ds i\kappa h})\,, \qquad
\hat y(\kappa,t_n;h) = \tilde y(e^{\,\ds i\kappa h},t_n) =
 \left(\tilde p(e^{\,\ds i\kappa h})\right)^n\,.
\leqno(3.9) $$
Recalling our fixation of
$t= t_n= n\tau = n\sigma(h) >0\,,$
the scaling relation (3.5) and the equivalences (3.6) we have
to show that
$$
\hat y(\kappa,t;h) \to \exp(-t|\kappa|^\alpha)\ \ \hbox{for} \ \ n\to
\infty\,,
\leqno(3.10) $$
or, equivalently
$$
{1\over \sigma(h)} \log \tilde p(e^{i\kappa h}) \to -|\kappa|^\alpha \ \
\hbox{as}\ \ h\to 0\,.
\leqno(3.11) $$
\vsp
In the following sections we shall exhibit (3.10) as true for
specific choices of the probabilities $p_j$ and scalings
$\tau = \sigma(h)\,.$ The fact that, strictly speaking,
$\tau = t/n= \sigma(h)$ and $h$ in
(3.6) and (3.11) are running through discrete sets will turn out as
irrelevant for the proof of (3.11).


\section{The Gr\"unwald-Letnikov random walk}
An idea suggesting itself is to discretize in (1.11) the time
derivative ${\partial u\over \partial t}$ by a two-level difference
quotient and the operators $I_+^{-\alpha}$ and $I_-^{-\alpha}$ (see
(2.15))
by the	Gr\"unwald-Letnikov approximation (see, \eg,  [SKM], [P]).
This idea leads to
$$
{y_j(t_{n+1}) - y_j(t_n)\over \tau} \,= \,_h D_0^\alpha y_j(t_n)
  = -{1\over 2\cos {(\alpha\pi/ 2)}}\,
 ( _h I_+^{-\alpha}+ _h I_-^{-\alpha}) \,y_j(t_n)\,.
\leqno(4.1) $$
with the operators $\,_h I_{\pm}^{-\alpha}$ still to be specified.
We must exclude the singular case $\alpha = 1\,,$ hence will distinguish
from now on the cases
$$\hbox{(a)}\q 0<\alpha<1\,,\qq \hbox{(b)} \q 1<\alpha\le 2\,.$$
\vsp
We define
$$
 _h I_\pm^{-\alpha} \,y_j (t_n) =
h^{-\alpha}\, \sum_{k= 0}^\infty (-1)^k
\l({\alpha \atop k}\r)\,y_{j\mp k}(t_n) \qq \hbox{in case (a)}\,,
\leqno(4.2) $$
$$
 _h I_\pm^{-\alpha} \,y_j (t_n) =
   h^{-\alpha}\,\sum_{k=  0}^\infty (-1)^k
\l({\alpha \atop k}\r)\,y_{j\pm 1 \mp k}(t_n) \qq \hbox{in case (b)} \,.
\leqno(4.3) $$
Note in case (b) the shift of the index, that is required in order
to obtain non-negative values for all transition probabilities $p_j\,.$
\vsp
Solving (4.1) for $y_j(t_{n+1})\,,$ thereby scaling by
$$
\tau =	   \mu h^\alpha =    : \sigma(h)\,,
\leqno(4.4) $$
gives (remember (3.4))
$$
y_j(t_{n+1}) =	\sum_{k\in\Z} p_k \,y_{j-k}(t_n)
\leqno(4.5) $$
with in case (a)
$$
p_0 =	 1-{\mu \over \cos{(\alpha \pi/ 2)}}\,,
\q
p_k=	(-1)^{|k|+1} {\mu \over 2\cos{(\alpha \pi / 2)}}\,
{\alpha \choose |k|}\q \hbox{for} \q k\not =0\,,
\leqno(4.6) $$
in case (b)
$$
\cases{
p_0  =	 1+{\ds{\mu \over \cos{(\alpha\pi/ 2)}}}\,
   {\ds{\alpha \choose 1}}\,,
  &\cr\cr
p_{\pm 1}  =
     -{\ds{\mu \over 2 \cos{(\alpha \pi/ 2)}}}\,
\left[ {\ds 1+ {\alpha\choose 2}}\right]\,,&\cr\cr
p_{\pm k}  =   (-1)^{k} \,
{\ds{\mu \over 2\cos{(\alpha \pi/ 2)}}}\,
{\ds{\alpha \choose k+1}} &
  for	$\;k=2,3, \dots\,.$\cr\cr}
\leqno(4.7) $$
Then,  all $p_k\ge 0$ if
$$
0<\mu \le \cos(\alpha\pi/2)\qq \hbox{in case (a)}   \,,
\leqno(4.8) $$
$$
0<\mu \le |\cos(\alpha\pi/2)|/\alpha\qq \hbox{in case (b)}  \,.
\leqno(4.9) $$
Note that $\cos(\alpha\pi/2)<0$ in case (b).
In both cases, by rearrangement of series,
$$
\sum_{k\in{\Z}} p_k = 1 - {\mu\over \cos{(\alpha\pi/ 2)}}\,
\sum_{j= 0}^{+\infty} (-1)^j
\,\l({\alpha \atop j}\r) = 1 - 0 = 1\,.
$$
{\bf Remark 4.1.}
For all $\alpha>0$ the series ${\ds \sum_{j= 0}^{+\infty}}\, (-1)^j\,
\l({\alpha \atop j}\r)$ is absolutely convergent because
$\l({\alpha \atop j}\r) =  O(j^{-\alpha -1})$
for  $j\to\infty\,.$
\vsp
We see that, under the conditions (4.8) or (4.9),
respectively, we can put
$$
P(Y= k) =     p_k\qquad (k\in \Z)
$$
for the random variable $Y$ of (3.1).
Using (4.6) and (4.7) we identify
the generating function $\tilde p$ of (3.7)  as
$$
\tilde p (z)
=     1 - {\mu\over 2\cos{(\alpha\pi/2)}}
\big\{(1 - z)^\alpha + (1- z^{-1})^\alpha\big\}
\qq \hbox{in case (a)},
\leqno(4.10) $$
$$
\tilde p (z)
=  1 - {\mu\over 2\cos{(\alpha\pi/ 2)}}\, \big\{z^{-1}(1 - z)^\alpha +
z(1- z^{-1})^\alpha\big\}\q \hbox{in case (b)}.
\leqno(4.11) $$
\vsp
Let us verify the limit relation (3.11) which implies (3.10).
Because of the symmetry relation (for $z= e^{i\kappa h},\ \kappa\in\R$)
$$
\tilde p(z) =	\tilde p(z^{-1})\ \ \hbox{implying}\ \
\tilde p(e^{i\kappa\,h})= \tilde p( e^{-i\kappa h})
$$
it suffices to verify (3.11) for $\kappa >0$ (the special case
$\kappa =    0$ being trivial).
\vsp
Let be $\kappa >0$. Then in case (a) we have
$$
\tilde p(z) =  1 -{\mu\over \cos{(\alpha\pi/ 2)}}\,\Re (1-z)^\alpha
  = 1-\mu (\kappa h)^\alpha +o(h^\alpha)  \q \hbox{as}\q h\to 0\,,
$$
since
$(1-z)^\alpha \sim (-i\kappa h)^\alpha =
  e^{-i\alpha\pi/2}\,(\kappa\,h)^\alpha\,.$
With the scaling (4.4), namely
$\sigma(h)= \mu h^\alpha\,,$ follows (3.11).
\vsp
In case (b) we have
$$
\tilde p(z) =	  1 -{\mu\over \cos{(\alpha\pi/ 2)}}\,
\Re (z^{-1}(1-z)^\alpha)\,,
$$
and an analogous  calculation gives
$$
{1\over \sigma(h)}\log \tilde p(e^{i\kappa h}) \sim
-\kappa^\alpha {\cos[{(\alpha\pi/ 2)}+\kappa h]\over
\cos{(\alpha\pi/ 2)}}\q \hbox{as}\q h\to 0\,,
$$
hence again (3.11). As result we have
\vsp  
{\bf Theorem 4.2.}
{\it Distinguish the cases
${(a)}\q  0<\alpha<1\,,
\q { (b)} \q 1<\alpha\le 2\,.$
Define the probabilities $p_k = P(Y= k)$ in case {(a)}
 by {(4.6)} with restriction {(4.8)},
in case {(b)} by
{ (4.7)} with restriction {(4.9)}.
Let the scaling relation $\tau= \mu \,h^\alpha = \sigma(h)$
hold and let for
fixed $t>0$ the index $n= t/\tau$ run through $\N$ towards $\infty$.
Then
the random variable $S_n$ of { (3.3)} converges in distribution
to the random variable $S(t)$ whose probability density is given by
{ (1.15)} as $g_\alpha(x,t;0) \,. $
} 
\vsp
{\bf Remark 4.3.}
In the special case $\alpha = 2$ the familiar explicit difference scheme
$$
y_j(t_{n+1}) =(1-2\mu) \,y_j(t_n) + \mu \,y_j(t_{n-1}) +
   \mu \, y_j(t_{n+1})
$$
is recovered from \hbox {(4.7)}, and \hbox{(4.9)} goes over
into the well-known stability condition $0<\mu \le 1/2\,.$
\vsp 
{\bf Remark 4.4.}
The case $\alpha = 1$ is singular.
For $\alpha \to 1$ both upper bounds in
\hbox{(4.8)} and \hbox{(4.9)} tend to 0, and the
denominators
occurring in the definitions of the probabilities $p_k$ tend to zero.
\vsp 
{\bf Remark 4.5.}
A motivation for the Gr\"unwald-Letnikov approximation of $I_+^{-\alpha}$
can be drawn from the fact that $z= e^{i\kappa h}$ is the symbol of the
backward shift by a step $h$: With
$$
(T_h\phi )(x) =    \phi(x+h)\,,\qq (T_{-h}\phi )(x)
=    \phi(x-h)\,,
$$
we have
$$
\widehat{T_h}(\kappa) =    e^{-i\kappa(x+h)}e^{i\kappa x}=
e^{-i\kappa h}\,,\q
\widehat{T_{-h}}(\kappa) =     e^{i\kappa h}\,.
$$
From the symbol
$\,h^{-1}(1-z) = \,_h \widehat{D_+} (\kappa)\,$
of the usual backward approximation
$$
(\,_h  D_+ \phi)(x) =
h^{-1}(\phi(x) -\phi(x-h)) =  h^{-1}\,(_h I \phi)(x)
$$
we arrive by analogy at the
symbol $h^{-\alpha}(1-z)^\alpha$
as a candidate for the symbol of the
operator $\,_h I_+^{-\alpha}\,.$ Analogously
we get $h^{-\alpha}(1-z^{-1})^\alpha$ as the symbol
for the operator $\,_h I_-^{-\alpha}\,.$
We use	the corresponding
approximations in  case \hbox{(a)} $\,0<\alpha<1\,.$
\vsp
In case \hbox{(b)} $\,1<\alpha\le 2$ we use the form
$D^2\,I_+^{2-\alpha}$
of the Riemann-Liouville left inverse of the operator $I_+^\alpha$,
and put $\,_h I_+^{-\alpha}\,=\, _hD^2\,_h I_+^{2-\alpha}\,.$
The corresponding symbol then is, with symmetrically
$(\,_h D^2 \phi)(x) = h^{-2}\,(\phi(x+h)-2\phi(x)+\phi(x-h))\,,$
$$
\,_h \widehat{I_+^{-\alpha}}(\kappa)\,=\,_h
\widehat{D^2}(\kappa)
 \,_h \widehat{I_+^{2-\alpha}}(\kappa) =
h^{-2}\,(z^{-1}- 2+z)\,h^{2-\alpha}\,(1-z)^{-(2-\alpha)}
$$
$$
=    h^{-\alpha} \,z^{-1}\,(1-2z+z^2)\,(1-z)^{\alpha -2} =
h^{-\alpha}\,z^{-1}\,(1-z)^\alpha\,.
$$
\vsp
The symbol $h^{2-\alpha}\,(1-z)^{\alpha -2}$ for
$\,_h I_+^{2-\alpha}$ here has been derived by the formal stipulation
$\,_h I_+^{2-\alpha}\,= \,_h D_+^{-(2-\alpha)}$ using
$\,_h \widehat{D_+} (\kappa) =	h^{-1}\,(1-z)\,.$
\vsp
Analogously we get
$ h^{-\alpha}\,z\,(1-z^{-1})^\alpha$
as symbol of $\,_h \widehat{I_-^{-\alpha}}\,.$
\vsp
{\bf Remark 4.6.}
In [GM98] and [GM99] we have exploited the
Gr\"unwald-Letnikov random walks in the more general setting of not
necessarily symmetric L\'evy-Feller diffusion (see Section 1).
The  proof for the case of symmetry ($\theta=0$)
given in the present paper is considerably simpler.


\section{The Gillis-Weiss random walk}
Gillis and Weiss in 1970 (see [GiW]) showed
(we interpret one of their results in the language of probability theory)
that every symmetric random variable $Y$
with values in $\Z$ and asymptotically
$P(Y=k)\sim {c/ |k|^{\alpha+1}}$  (where $c>0$)
lies in the domain of attraction of the corresponding symmetric
L\'evy distribution, hence can be used for an approximating random walk
in the sense of Section 3.
Only assuming their asymptotics they naturally cannot describe
precisely how the coefficients
$\mu$ and $\lambda$ of the scaling law appear in the
transition probabilities.
However, from their
analysis we can deduce that the scaling law is of the form
$$
\tau =\sigma(h)= \mu \,h^\alpha\ \ \hbox{if}\ \ 0<\alpha<2\,,
 \qq
\tau =\sigma(h)=\lambda \, h^2\, |\log h|\ \ \hbox{if}\ \ \alpha=2\,.
$$
Remarkably, the parameter value $\alpha = 1$ is not singular, but the
scaling law becomes discontinuous at $\alpha =2$,
thus giving an example of a distribution with non-finite variance
lying in the domain of
attraction of the normal (Gauss) distribution.
\vsp
We will now re-work and complement their analysis in the
framework of our Section 3 for the special symmetric
probability distribution
$(p_k|\,k\in \Z$) with
$$
p_0=1-2\lambda \sum_{k=1}^\infty k^{-(\alpha +1)}\,,
\q p_k=    \lambda
|k|^{-(\alpha +1)}\q \hbox{for}\q k\not = 0\,,
\leqno(5.1) $$
where (so that $p_0 \ge 0$) $\lambda$ is restricted by
$$
0<\lambda \le
{\left(2\,{\ds \sum_{k=1}^\infty  k^{-(\alpha+1)}}\right)}^{-1}\,.
   \leqno(5.2) $$
The parameter $\alpha$ is only restricted as in (1.4) by
$0<\alpha\le 2\,.$ Differently from Gillis and Weiss we motivate this
choice of probabilities by (2.20), where the special character of
the value $\alpha = 2$ already becomes visible. So, assume meanwhile
$0<\alpha<2\,.$
\vsp
Discretizing $D_0^\alpha u$ via a straightforward quadrature
formula for the right hand side of (2.20) as
$$
\,_h D_0^\alpha y_j(t_n) =  b(\alpha) \,h \,\sum_{k=	1}^\infty
{y_{j+k}(t_n)-2y_j(t_n) +y_{j-k}(t_n)\over (kh)^{\alpha +1}}
\leqno(5.3) $$
and solving the equation
$$
{y_j(t_{n+1}) - y_j(t_n)\over \tau} \,=\,_h D_0^\alpha\, y_j(t_n) $$
for $y_j(t_{n+1})$
we identify the transition probabilities $p_k$ in (3.4) as
$$
p_0=   1-2\mu \,b(\alpha) \,\zeta(\alpha +1),\
p_k=	\mu b(\alpha)\, |k|^{-(\alpha +1)}\ \ \hbox{for}\ \ k\not = 0\,,
\leqno(5.4) $$
with $\mu =  h^{-\alpha}\tau\,, $
$ b(\alpha ) = \Gamma(\alpha+1)\, \sin {(\alpha \pi/2)}/\pi\,$
and the Riemann $\zeta$-function
$$
\zeta(z) = \sum_{k= 1}^\infty k^{-z}\ \ \hbox{for}\ \ \Re z >1.
\leqno(5.5) $$
Obviously $\sum_{k\in \Z} p_k = 1\,,$ and the non-negativity condition
in (3.2) requires
$$
0<\mu \le {1\over 2 \,b(\alpha)\,\zeta(\alpha +1)}=
{\pi\over 2\Gamma(\alpha+1)\,\sin{(\alpha \pi/2)}\,\zeta(\alpha+1)}
\,. \leqno(5.6) $$
\vsp
We want to free the parameter value $\alpha = 2$ from its singular
character. Recalling (2.17)
we see that $b(2)=0\,,$ so that in (5.4) $p_0=1$ and all $p_k=0$
for $k\not = 0$ whereas the upper bound for $\mu$ in (5.6) tends to
$\infty$ as $\alpha \to 2-$. This degenerate random walk obtained in
(5.4) by formally setting $\alpha = 2$ being neither interesting
nor useful  we replace $\mu\, b(\alpha)$ by $\lambda$ and obtain the
transition probabilities in the form (5.1) with restriction (5.2).
In (5.1) the special value $\alpha = 2$
seems to be a quite regular value, and we shall see that we
have a valid random walk model for all $\alpha$ obeying
$0<\alpha\le 2\,.$
However a price must be paid. Whereas for $0<\alpha <2$ we can scale by
$\tau = \mu \,h^\alpha$
we can no longer do so in the case $\alpha= 2\,.$
So, assume henceforth (if not explicitly stated otherwise) the condition
(5.2).
\vsp
We have now the generating function
$$
\tilde p(z) =  1-2\lambda \zeta(\alpha +1) +\lambda \sum_{k=1}^\infty
k^{-(\alpha+1)}(z^k + z^{-k})
\leqno(5.7) $$
with $z=  e^{i\kappa h},\ \kappa\in \R\,.$
With the polylogarithmic function
$$
\Phi(z,\beta) =  \sum_{k= 1}^\infty {z^k\over k^\beta}\qquad
 (\beta \in \R)
$$
we can write
$$
\tilde p(z) =  1-2\lambda \,\zeta(\alpha +1) +\lambda \,
\{\Phi(z,\alpha+1) +\Phi(z^{-1},\alpha+1)\}
\leqno(5.8)$$
and  could carry out the required asymptotic analysis by specializing
some of the formulas in [T]. See also [EHTF] and [Le] for properties of
the polylogarithmic function and the more general Lerch function. We
prefer, however, the direct way to obtain
(3.11). This asymptotic relation is trivial for $\kappa = 0\,,$ and
because of $\tilde p(e^{i\kappa h}) =  \tilde p(e^{-i\kappa h})\,,$
it suffices to treat the case $\kappa >0$ what we now will do.
\vsp
From the common integral representation of the gamma function we take
$$
k^{-(\alpha +1)} = {1\over \Gamma(\alpha +1)} \int_0^\infty u^\alpha
e^{-ku} du
$$
and using $z^{-1}= \bar z$ we get
$$
\tilde p(z) =	  1-2\lambda \,\Re \gamma(z)
\leqno(5.9) $$
with
$$
\Gamma(\alpha +1)\,\gamma(z) =
  \int_0^\infty u^\alpha \,\sum_{k= 1}^\infty e^{-ku}(1-z^k)\, du
=  \int_0^\infty {u^\alpha\, e^{-u} \over 1-e^{-u}}\,
{{1-z}\over 1-e^{-u}\,z}\,du\,.
\leqno(5.10)
$$
The last equality in (5.10) has been obtained
by  summing the two involved geometric series.
\vsp
In the Appendix we have performed in detail the required asymptotic
analysis of  $\Re \gamma(z)$ for $\nu = \kappa \,h \to 0+$ ($\kappa$
fixed), which is resumed in formulas (A.6) and (A.7).
Insertion of these asymptotic behaviours into (5.9) yields
$$
\log \tilde p(e^{i\kappa h})
\sim - {\lambda \pi \over \Gamma(\alpha +1)
\sin{(\alpha \pi/ 2)}}\, |\kappa|^\alpha h^\alpha \ \ \hbox{if}\ \
0<\alpha <2,\ \kappa\not = 0\,,
\leqno(5.11) $$
$$
\log \tilde p(e^{i\kappa h})\sim - \lambda \kappa^2 h^2
\, \log {\l(1/(|\kappa|h) \r)}\ \
\hbox{if}\ \ \alpha =	 2,\ \kappa\not =    0\,.
\leqno(5.12) $$
\vsp
Recalling that it suffices to prove (3.11)
for $\kappa\not = 0$ and
observing, that there the parameter $\kappa$ can be treated like a
constant, we see that $\log (1/(|\kappa|h))\sim \log{(1/ h)}\,,$
where, because $h\to 0\,,$ we can assume $0<h<1\,. $
\vsp
Hence we can replace (5.12) by
$$
\log \tilde p(e^{i\kappa h})\sim -\lambda \kappa^2 h^2 \log {1\over h}\ \
\hbox{if}\ \ \alpha =	 2,\ \kappa\not = 0\,.
\leqno(5.13) $$
Then the limit relation (3.11) (equivalently (3.10)) holds
if we scale by
$$
\tau =	\sigma(h) = {\lambda \pi \over \Gamma(\alpha +1)
\,\sin{(\alpha \pi/ 2)}} \, h^\alpha \q \hbox{if}\q
0<\alpha<2\,,
\leqno(5.14) $$
$$
\tau =	\sigma(h) =\lambda h^2 \log {1\over h} \q
 \hbox{if}\q \alpha =  2\,.
\leqno(5.15) $$
\vsp
Putting $\mu = \lambda \pi /(\Gamma(\alpha +1)\,
\sin{(\alpha \pi/2}) =	  {\lambda/b(\alpha)}$ in (5.9)
with $b(\alpha)$ defined in (2.17) we obtain from (5.4)
the regular scaling law
$$
\tau =	\sigma(h) = \mu \,h^\alpha\ \ \hbox{for}\ \ 0<\alpha <2
\leqno(5.16) $$
with the restriction  (5.6) for $\mu \,.$
As result we have
\vsp 
{\bf Theorem 5.1.}
{\it Distinguish the cases
{(i)} $\;0<\alpha<2\,,\q $ {(ii)} $\;\alpha =  2\,.$
Define the probabilities $p_k= P(Y= k)$ in case {(i)} by
{(5.4)} with restriction {(5.6)},
in case {(ii)} by 
$$
p_0 =	 1-2\lambda \zeta(3),\ \ p_k =	  \lambda |k|^{-3} \q
 \hbox{for} \q k \ne 0
$$
with restriction $0<\lambda \le {1/(2\zeta (3))}\,.$
Let the scaling relation
$$
\tau =	  \mu h^\alpha\ \
 {in} \ { case} \ {(i)},
\qq
\tau =	  \lambda h^\alpha \log{1\over h}\ \
{ in} \ { case} \  { (ii)}
\leqno(5.17)
$$
hold and let for fixed $t>0$ the index $n= t/\tau$ run through $\N$
towards $\infty\,.$ Then the random variable $S_n$ of {(3.3)}
converges in distribution to the random variable $S(t)$ whose probability
density is given by
{(1.15)} as $g_\alpha(x,t;0) \,. $}
\vsp
{\bf Remark 5.2.}
We can use throughout $0<\lambda \le 2$ the parameter $\lambda$ and then
have in  (5.1) under the restriction (5.2)
a unified representation of the transition probabilities.
Here, in
contrast to the Gr\"unwald-Letnikov random walk, the value $\alpha =1$
does no longer play a special role. With $\mu = \lambda/b(\alpha)$
we have for $0<\alpha <2$ the regular scaling law $\tau = \mu h^\alpha\,.$
However, the price to be paid for this unified representation
is the non-regular scaling
$\tau = h^2\,\log (1/h)$ for $\alpha= 2\,.$
Another price is that the generating function $\tilde p(z)$ in (5.8)
is non-elementary, requiring considerable efforts in its asymptotic
analysis.


\section{A globally binomial random walk}
The random walk model discussed in Section 4 has the disadvantage that
the case $\alpha =1$ is excluded and the representation of the
transition probabilities $p_k$ for $1<\alpha \le 2$ is different from that
for $0<\alpha <1\,. $
However, for all admissible values of $\alpha$ we have the
regular scaling law $\tau =    \mu h^\alpha$.
The method treated in Section 5
has the advantage of a unified representation of the transition
probabilities in the whole interval $0<\alpha \le 2\,,$
but the scaling law $\tau=    \mu \,h^\alpha$
holds only for $0<\alpha <2\,, $
it breaks down at $\alpha =  2\,.$
In this section we present a model that in the whole
interval $0<\alpha \le 2$ admits a unified representation of the $p_k$
via binomial coefficients
and has there a scaling law of the form $\tau= \mu \,h^\alpha\,.$
Moreover,
the generating function $\tilde p (z)$
is elementary for all $\alpha \in (0,2]\,.$
\vsp
The use of the binomial coefficients ${\alpha \choose j}$ in the
Gr\"unwald-Letnikov random walk has caused singular behaviour
for $\alpha=1\,.$
 One reason for this sad fact is that
${1\choose j} = 0$ for integer $j\ge 2\,.$
We can remove this singular behaviour by removing the factor
$\alpha -1$.
\vsp
For $0<\alpha\le 2\,,\ \alpha \not = 1\,$ let us define
$$
p_0=	1-2\lambda,\ p_k =
 (-1)^{k+1} {\lambda\over \alpha -1} \,{\alpha \choose |k|+1}\ \
\hbox{for}\ \ k\not =	 0\,.
 \leqno(6.1) $$
Observing that here the singularity at $\alpha = 1$ is
removable, let us for $\alpha = 1$ define
(via $\alpha \to 1$ in (6.1))
$$
p_0=	1-2\lambda\,,\q p_k =	 {\lambda\over |k|(|k|+1)}\ \
\hbox{for}\ \ k\not =	0\,.
\leqno(6.2) $$
In (6.1) and (6.2) $\sum_{k\in\Z}p_k = 1$ and if $0<\lambda\le 1/2$
all $p_k\ge 0\,.$ In the special case $\alpha = 2$ we get
$$
p_0= 1-2\lambda\,,\q p_1=  p_{-1}= \lambda\,,\q p_k=	0 \ \
\hbox{for}\ \ |k|\ge 2\,,
$$
the familiar random walk for approximation of the classical process
governed by the equation
${\partial u\over \partial t}= {\partial^2 u\over \partial x^2}\,.$
\vsp
The generating function $\tilde p(z)= \sum_{k\in \Z}p_k z^k$ has in the
case $\alpha \not =    1$ the form
$$
\tilde p(z) =	 1-\lambda \{ q(z)+q(z^{-1})\}
\leqno(6.3) $$
with
$$
q(z) =	{1\over \alpha -1}\,(1-z^{-1})\,\{ (1-z)^{\alpha -1}-1\}\,.
$$
By passing here to the limit or directly from (6.2) we get for
$\alpha =    1$ the representation
$$
\tilde p(z) =  1 -\lambda \{(1-z^{-1})\log (1-z) +(1-z)\log(1-z^{-1})\}
 \,, \q \tilde p(1) = 1\,.
\leqno(6.4) $$
We have proposed and investigated the particular random walk
so generated (its transition probabilities given in (6.2))
in [GM99, Section 5].
\vsp
In the special case $\alpha =	 2$ we find
$$
\tilde p(z) =	1 +\lambda (z-2+z^{-1})\,.
\leqno(6.5) $$
\vsp
We will now show that for all $\alpha \in (0,2]$ there
exists a finite positive number $c(\alpha)$ so that, with
$$
\mu = c(\alpha)\,\lambda   \,,
\leqno(6.6) $$
we arrive for  $\kappa \in \R\setminus \{0\}$ at the small $h$
asymptotics
$$
\tilde p(e^{i\kappa h}) = 1-\mu (|\kappa| h)^\alpha +
   o\left( (|\kappa|\, h)^\alpha \right)
\leqno(6.7) $$
which implies (3.11). As in Sections 4 and 5 we can ignore the
value $\kappa =    0$ as trivial.
\vsp
Referring to [GM99] for detailed treatment of the case $\alpha =1\,, $
let now be $0\not =  \kappa \in \R$ and
$ 0<\alpha \le 2,\ \alpha \not =  1\,,\, z=  e^{i\kappa h}$.
In view of (6.3) we investigate the asymptotics
of $q(z)+q(z^{-1})$ for $h\to 0$. From
$z^{-1} = \bar z$ and
$$
(1-\alpha)q(z) =  z^{-1}(1-z)^\alpha - z^{-1} +1 =
    e^{-i\kappa\,h}(1-e^{i\kappa\, h})^\alpha - e^{-i\kappa\, h} +1\,,
$$
we conclude on
$$
\psi(z):= (1-\alpha)\{q(z)+q(z^{-1})\} =
 2\Re\left\{e^{-i\kappa h}(1-e^{i\kappa h})^\alpha\right\} +
  2(1-\cos(\kappa h))\,,  \leqno(6.8) $$
and here
$$
\Re\left\{e^{-i\kappa \,h}(1-e^{i\kappa\, h})^\alpha\right\}
\sim \Re\left((-i\kappa\, h)^\alpha\right) =
   (|\kappa|\,h)^\alpha \cos{(\alpha\pi/2)},
\leqno(6.9) $$
$$
1-\cos(\kappa h)\sim {1\over 2}(|\kappa |h)^2\,.
\leqno(6.10) $$
\vsp
We distinguish three cases:
$
\hbox{(i)}\ \ 0<\alpha<1\,,
\ \ \hbox{(ii)}\ \ 1<\alpha<2\,,\ \
\hbox{(iii)} \ \alpha =    2\,.
$
\vsp
In cases (i) and (ii) the leading term
in the asymptotics of $\psi(z)$ turns out to be
$$\psi(z)\sim 2\,(|\kappa|\,h)^\alpha\,\cos{(\alpha\pi/ 2)}\,.$$
\vsp
In case (iii) where $\alpha = 2$ however, this term is matched in order of
magnitude by (6.10) so that we obtain
$$
\psi(z)\sim 2(|\kappa|\,h)^2(-1)+(|\kappa|\,h)^2= -(|\kappa|\,h)^2\,.
$$
\vsp
Collecting results and dividing (6.8) by $1-\alpha$ we get
(with $z= e^{i\kappa h}$)
$$
\lambda\{ q(z)+q(z^{-1})\}\sim
\cases{
\lambda {\ds{2 \cos{(\alpha\pi/ 2)}\over 1-\alpha}}\,
(|\kappa|\,h)^\alpha & if
$0<\alpha<2\,,\ \alpha \not = 1\,,$\cr\cr
\lambda (|\kappa|\,h)^2 & if $\alpha =	2\,.$\cr}
\leqno(6.11) $$
Hence, in view of (6.3), we obtain (6.7)
with (6.6) by putting
$$
c(\alpha)=  \cases{
{\ds {2 \cos{(\alpha\pi/2)}\over 1-\alpha}}
   & if $\q 0<\alpha<2,\  \alpha \not = 1\,,$\cr\cr
 1 & if $\q \alpha = 2\,.$\cr}
\leqno(6.12) $$
\vsp
The scaling coefficient $c(\alpha)$ allows continuous extension to the
value $\alpha = 1\,,$ giving
$\lim\limits_{\alpha \to 1} c(\alpha) = \pi$
in accordance with [GM99, formula (5.1)].
At $\alpha = 2\,,$ however, $c(\alpha)$ is discontinuous. In fact
$$
c(2)=	 1\not =   2= \lim_{\alpha \to 2} c(\alpha)\,.
\leqno(6.13) $$
\vsp
Let us finally display the transition probabilities with $\mu$ instead of
$\lambda$ as parameter.
\vsp
For $0<\alpha<2,\ \alpha \not = 1$:
$$
\cases{p_0 = 1-2\mu \,
{\ds {1-\alpha\over 2\,\cos{(\alpha\pi/ 2)}}}\,,&\cr\cr
p_k={\ds{(-1)^k\over 2\cos{(\alpha\pi/2)}}}\,{\ds{\alpha\choose |k|+1}} &
for   $\;k\not= 0\,,$\cr\cr
0<\mu\le {\ds{\cos{(\alpha\pi/2)}\over 1-\alpha}}\,,& \cr }
\leqno(6.14) $$
for $\alpha = 1$ (see [GM99, formula (5.1)]):
$$
 p_0 = 1-{\ds{2\mu \over \pi}}\,,\q
 p_k= {\ds {\mu\over \pi |k|(|k|+1)}} \q \hbox{for}  \q k\not = 0\,,\q
 0<\mu\le {\pi/ 2}\,,
\leqno(6.15) $$
for $\alpha =  2$:
$$
p_0 =  1-2\mu\,,\q
p_1= p_{-1}= \mu,\q
p_k= 0 \q \hbox{for} \q |k|\ge 2\,,\q
0<\mu \le 1/2\,.
\leqno(6.16) $$
The discontinuity at $\alpha = 2$ has so been transferred to the upper
bound for $\mu\,.$
\vsp
We comprise the result in
\vsp
{\bf Theorem 6.1.}
{\it Take the probabilities $p_k=  P(Y= k)$ and the	restrictions for $\mu$ as
in formulas
{(6.14)}, { (6.15)}, {(6.16)}, and use
the scaling relation
$\tau =  \mu \,h^\alpha\,.$
Let for fixed $t>0$ the index $n=t/\tau$ run through $\N$ towards
$\infty\,. $ Then the random variable $S_n$ of { (3.3)} converges
in distribution to the random variable $S(t)$ whose probability
density is given by {(1.15)} as $g_\alpha(x,t;0) \,. $}


\section{The Chechkin-Gonchar random walk}
In this section we adopt to each other considerations of Chechkin
and Gonchar [ChG] and the framework of our Section 3, restricting attention
to the parameter range $ 0 <\alpha <2\,.$
So doing
we exclude the well-known case of the classical Gaussian process.
We will obtain a {\it random walk}, which is {\it discrete in time}
but {\it continuous in space}, in more precise words:
whose jumping width (in the instants $t_n= n \tau$) can assume any
real number, having an everywhere positive probability density.
We modify our theory of Section 3 by allowing the random variable
$Y$ to have a strictly monotonic continuous distribution function
$ W(x) = P\l( Y<x \r)\,$ $\,(x \in \R),$
that we furthermore require to be symmetric in the sense
$ W(x) + W(-x) = 1\,$ $\,(x \in \R),$
being only interested in the symmetric case $\theta=0$ of Section 1.
We then have in the sum
$$
S_n= \sum_{m=1}^n \l(h\,Y_m \r)= h\,\sum_{m= 1}^n Y_m \qq (n\in \N)
\leqno(7.1)$$
a description of a random walk, starting in the point $x=0\,.$
Here $h>0$ is a scaling width that we let depend on the time-step
$\tau >0$ via a strictly monotonic scaling relation
$\tau =\sigma (h)\,, $ with $\sigma (h) \to 0$ as $h \to 0\,. $
We expect  the scaling relation  to have the form $\tau =\mu \,h^\alpha $
with the positive coefficient $\mu $ to be specified,
by having found orientation in Gnedenko's theorem on
normal attraction (see [GnK], \S 35).
It should be noted, however,  that in this theorem the scaling constant
$C$ appearing there is given with a wrong value as has been
remarked in [Ba].
\vsp
As previously, we let  the $Y_m$ be independent identically distributed
random variables,
all having their distribution common with $Y\,.$
However, we now assume $Y$ to have an everywhere positive
(not necessarily bounded)  probability density
$w = W'$ which  is an even function
$ w(x) = w(-x)\,$ $\, (x\in\R)\,. $
We will use the fact that $w$ is normalized,
$ \int_{-\infty}^{+\infty} w(x)\, dx =1\,.$
\vsp
Fixing a value $t>0$ and again setting $t=t_n = n\tau$
(equivalent to $n=t/\tau $) with $n\in \N$ we want that the
random variable $S_n$ converges in distribution to the
random variable $S(t)$ whose density is given by (1.15).
To this purpose we introduce a condition on the asymptotic behaviour
of the density $w$, namely
$$ w(x) = \l( b+ \epsilon (|x|)\r) \, |x|^{-(\alpha +1)}\,,
\q   |\epsilon (|x|)| \le \hbox{min} \,\l\{K, E\,|x|^{-\gamma}\r\}
 \q (x\in \R),\leqno(7.2)$$
with positive constants $b\,,\,K\,,\,E$ and $\gamma \,. $
\vsp
With $\hat w(\kappa) = \int_{-\infty}^{+\infty} e^{i\kappa x}\,w(x)\,dx$
as characteristic function of the density $w(x)$ we observe that
the random variable $hY$ has density $w(x/h)/h\,,$ hence the characteristic
function $\hat w(\kappa h)\,,$ and proceeding in analogy to the general
method described in Section 3, replacing
$\hat p(\kappa ,h) = \tilde p (e^{i\kappa h})$ in (3.9)
by $\hat w(\kappa h)\,, $ we will find a scaling function
$\sigma (h)$ such that for all $\kappa \in\R\,,$
in analogy to (3.11),
$$ {1\over \sigma (h)}\, \log \hat w(\kappa h) \to
  - |\kappa |^\alpha \q \hbox{as}\q h \to 0\,. \leqno  (7.3)$$
Of course, (7.3) is trivial for $\kappa =0\,. $
Since $\hat w$ like $w$ is an even function
it suffices to consider (7.3) for (fixed) values $\kappa >0\,. $
 We see that (7.3) is equivalent to
$$ \hat w(\kappa h) = 1- |\kappa |^\alpha \, \sigma (h) + o(\sigma (h))
  \q \hbox{as}\q h\to 0\,. \leqno(7.4)$$
\vsp
In view of the symmetry and normalization properties of $w(x)$
and abbreviating $\kappa h=\nu $ we find
$$\hat w(\nu) -1
= \int_0^\infty \l(e^{i \nu x} + e^{-i \nu x} -2\r) \, w(x)\, dx
= -4 \,  \int_0^\infty \l(\sin(\nu x/2)\r)^2 \,w(x)\, dx$$
so, using (7.2),
$$ \hat w(\nu) = 1
-2^{-\alpha +2} \, b\, \nu^\alpha \,\int_0^\infty
  \xi ^{-\alpha -1} \, \l(\sin \xi \r)^2 \, d\xi
 -4 \int_0^\infty \epsilon (x) \,  x^{-\alpha -1}\,
  \l(\sin(\nu x/2)\r)^2 \, dx\,.
 $$
The first integral can be evaluated in terms of the gamma function.
In fact, from [GR, (3.823)] we take
$$ \int_0^\infty  \xi^{-\alpha -1} \, \l(\sin \xi \r)^2 \, d\xi
= - {\Gamma(-\alpha )\,\cos (\alpha \pi/2)\over 2^{1-\alpha}}
 = {\pi \over 2^{2-\alpha }\,\Gamma (\alpha +1)\, \sin(\alpha \pi/2)}
\,. $$
The latter equality follows by the reflection formula
for the gamma function.
\vsp
We estimate the  second  integral
via decomposition
$\int_0^\infty \dots = \int_0^\eta \dots  + \int_\eta^\infty \dots\,,$
taking $\eta = \nu ^{-(2\alpha +\gamma) /(2\alpha +2\gamma )}$,
using	$|\sin \xi| \le \hbox{min}\,\{\xi ,1\}$ for $\xi \ge0$
and the condition on $\epsilon (|x|)$ of (7.2).
By careful calculation we  find that it behaves asymptotically
as $  o\l(\nu^\alpha\r) = |\kappa|^\alpha \,o\l(h ^\alpha\r)\,.$
Combining these results and recalling that $\hat w$ is an even function,
we obtain
$$ \hat w(\kappa h) = 1- |\kappa |^\alpha \,
   {b\,\pi \over \Gamma(\alpha +1)\,\sin(\alpha \pi/2)}\, h^\alpha
    + |\kappa |^\alpha \,o\l(h^\alpha \r)
\q (h\to 0)\leqno(7.5)$$
as valid for all $\kappa \in \R\,. $
In view of (7.4) and the theory developed in Section 3
we thus arrive at the scaling relation
$$ \tau = \sigma(h) = \mu \, h^\alpha \,,\q
   \hbox{with} \q \mu =
   {b\,\pi \over \Gamma(\alpha +1)\,\sin(\alpha \pi/2)} \,.
  \leqno(7.6)$$
Now we are in the position to formulate
\vsp 
{\bf Theorem 7.1.}
  {\it Let $0<\alpha <2$ and assume the random variable $Y$ to have a
probability density $w$ of the form { (7.2)}.
Let the scaling relation {(7.6)}
hold and let for fixed $t>0$ the index $n=t/\tau $ run
through $\N$ towards $\infty\,.$
Then the random variable $S_n$ of { (7.1)} converges in
distribution to the random variable $S(t)$ whose
probability density is given by {(1.15)}.}
\vsp
{\bf Remark 7.2.}
According to the well-known asymptotic expansions
of the function $p_\alpha (x;0) = g_\alpha (x,1;0)$
(see [F52], [F71], [Zo]) we have
$$ b = {\Gamma(\alpha +1)\, \sin(\alpha \pi/2)\over \pi}
\leqno(7.7)$$
if we take $w(x) = p_\alpha (x;0\,, $ hence in this case
$t=1\,,$ $\mu =1$ and $h=\tau ^{1/\alpha}= n^{-1/\alpha}\,. $
If we require in (7.1) the $Y_m$ to have this special
density, then $S_n$ for all $t>0$ has the same probability distribution
as $S(t)$ whose characteristic function is
$\exp (-t\,|\kappa|^\alpha )\,. $
We can here obtain the scaling relation also
via the  convolution theorem.
\vsp
{\bf Remark 7.3.}
For  actual simulation	a random variable $Y$ having the required
properties is particularly useful if
its distribution function
$W(x) = \int_{-\infty}^x  w(\xi )\,d\xi $
is easily invertible.
We can then generate a realization of $Y$ by a standard
Monte Carlo method (see [HH]). Generate a random number
$y$  uniformly distributed  in the  interval $[0,1)\,. $
Then solve the equation $y=W(x)$ for $x$ and take $x$
as a realization of $Y$. Chechkin and Gonchar in [ChG]
have proposed to use
$$ W(x) = \cases{
{1\over 2} \,\l(1+|x|^\alpha \r)^{-1}  & for $\;x<0\,, $\cr\cr
1-{1\over 2} \,\l(1+x^\alpha \r)^{-1}  & for $\;x\ge 0\,, $\cr}
\leqno(7.8) $$
a function easily invertible.
The density
$$ w(x) = W'(x) = {\alpha |x|^{\alpha -1} \over
 2 \,\l(1+|x|^\alpha \r)^2}\leqno (7.9)$$
has the property (7.2) with $b=\alpha /2\,,$ $\, \gamma =\alpha \,, $
hence we get
$$ \mu = {\pi \over 2\Gamma(\alpha)\, \sin(\alpha \pi/2)}\,.
\leqno(7.10)$$
The density (7.9) is unbounded at the origin if $0<\alpha <1\,. $
To avoid this we propose
 $$ w(x)  = {\alpha \over 2}\, \l(1+ |x|\r)^{-(\alpha +1)}
 \,, \leqno (7.11)$$
which again satisfies the asymptotic condition	(7.2).
Then
 $$ W(x) = \cases{
{1\over 2} \,\l(1+|x| \r)^{-\alpha}  & for $\;x<0\,, $\cr\cr
1-{1\over 2} \,\l(1+x\r)^{-\alpha }  & for $\;x\ge 0\, $\cr}
\leqno(7.12) $$
is also easily invertible, and (7.10) remains valid.
\vsp
{\bf Remark 7.4.}
Among the symmetric densities $p_\alpha (x;0)\,,$ only the
Cauchy density
$ w(x) = p_1(x;0) = ({1/ \pi})\, (1+x^2)^{-1}\, $
offers easy invertibility of the corresponding distribution
function, namely of the function
$ W(x)= {1/ 2} + ({1 / \pi})\, \hbox{arctan} \,x\,.$
Via random numbers $y_m$  uniformly distributed  in [0,1)
we can get  realizations of the $Y_m$ in
$\tau \, \hbox{tan} (\pi(y_m-1/2))$ (here $\tau = 1\, h^1 = h$)
and so obtain in $S_n$
a snapshot at instant $t_n = n\,\tau \,$ of a true
Cauchy process.
\vsp
\section{Conclusions}
Anomalous diffusion processes have in recent years gained revived interest
among physicists, and methods of fractional calculus have shown their
usefulness for purposes of modelling.
In the space-fractional case one is naturally led to a
generalization of the classical diffusion equation with respect
to the second-order spatial operator.
One arrives in a natural way at the
processes of L\'evy-Feller type in which stable probability
distributions play the essential role.
Also among physicists and mathematicians who have found it
rewarding to work in theory of finance,
such processes are becoming more and more
popular (see e.g.  [M], [BoP], [MS]).
So, it is no wonder that also in pure mathematics such types of
processes are now investigated in great generality
and analytical sophistication (see e.g. [J], [Be],  [S], [Za]).
From the more practical point of view discrete models are esteemed.
They not only show that very different microscopic behaviour
of particles can result in the same macroscopic behaviour but
offer also possible visualizations of what is happening in
such processes.
Furthermore such discrete models can be used for simulation purposes,
be it for simulation of particle paths via Monte Carlo
methods (the microscopic view) or via solution of the underlying
Cauchy problem for a pseudo-differential equation
(the macroscopic view). And, last but not least,
such models are fascinating as seen from the
mathematical standpoint (or, more specifically, from
the position of probability theory).
\vsp
In our present investigation we first have given a survey on and
drawn motivations from basic theory of fractional calculus and
L\'evy-Feller  diffusion processes.
Then we have obtained and rigorously analyzed
(with respect to their convergence in distribution for passing
to the limit of infinitely fine discretization)
three models of random walk occurring on a regular spatio-temporal grid.
The first model is devised   from the Gr\"unwald-Letnikov
discretization of the two Weyl operators, the composition of
which gives the inverse of the Riesz potential operator.
The second model is an adaptation of ideas of Gillis and Weiss [GiW]
to our framework. We have provided it with a new motivation,
namely as obtainable from straightforward discretization
of the hypersingular integral representation of the spatial
pseudo-differential operator.
The third model's intention is to overcome peculiar deficiencies of
the first two models. It is  a modification and improvement
of the first model, and again properties of the binomial coefficients
are used.
\vsp
Finally, to offer also a highly efficient method for numerical
simulation, we have mutually adapted our theoretical frame
to ideas of Chechkin and Gonchar [ChG]. We so obtain a random walk
still proceeding in equidistant instants of time but allowing
spatial jumps of arbitrary length in positive or negative
direction.
\vsp

\noindent
\section*{Appendix A: Asymptotics of an integral}
\vsp
Abbreviating $\kappa h= \nu$ in $z = e^{i\kappa h}$
in the right hand side of (5.10),
and keeping in mind $0 <\alpha \le 2\,, $
elementary calculation	yields the equation
$$
\Gamma(\alpha+1)\,\Re \gamma(z) =  \int_0^\infty u^\alpha \,e^{-u}\,
{(1+e^{-u}) \,(1-\cos \nu) \over
   (1-e^{-u})\, |1-e^{-u}\,e^{i\nu}|^2}\,du
\leqno(A.1) $$
which we will treat asymptotically for $0< \nu \to 0+$ by the Laplace
method
for integrals (see [dB]), using the fact that the lower bound $u= 0$ is
the critical one (the integrand tending to $\infty$ as $u\to 0$).
 We have $1-\cos \nu = {\nu^2/ 2}+O(\nu^4)$ and
$$
|1-e^{-u}e^{i\nu}|^2 = (1-e^{-u})^2 +2e^{-u}(1-\cos \nu) =
    (1-e^{-u})^2 +\nu^2e^{-u} +O(\nu^4),
$$
uniformly in $0\le u<\infty$, hence
$$
\Gamma(\alpha+1)\,\Re \gamma(z) \sim {\nu^2\over 2}
\int_0^\infty {u^\alpha e^{-u}(1+e^{-u})\over (1-e^{-u})\{
(1-e^{-u})^2 +\nu^2 e^{-u}\} }\, du\,.
$$
Because this integral diverges for $\nu = 0$ we can simplify the integrand
(for {\it small} $u$) which, for small $\nu$, gives the essential
contribution:
$1+e^{-u}\sim 2\,,$
$\, 1-e^{-u}\sim u\,,$
$\, e^{-u}\sim 1\,.$
We obtain
$$
\Gamma(\alpha+1)\,\Re \gamma(z) \sim {\nu^2}
\int_0^\infty u^{\alpha-1} {e^{-u}\over u^2+\nu^2}\,du
$$
and, by substituting $u= \nu \,w$,
$$
\Gamma(\alpha+1)\Re \gamma(z) =     \nu^\alpha \int_0^\infty
w^{\alpha-1} {e^{-\nu w}\over w^2+1}\,dw \, :=
   \, \nu^\alpha \rho(\nu)\,.
\leqno(A.2) $$
In the investigation of the integral
$$
\rho(\nu) =	\int_0^\infty
w^{\alpha-1}\, {e^{-\nu w}\over w^2+1}\,dw  \leqno(A.3)
$$
we distinguish the cases
$
\hbox{(i) }\ 0<\alpha<2,\quad \hbox{(ii) }\ \alpha = 2\,.
$
In the case (i) simply
$$
\rho(\nu) \rightarrow \int_0^\infty {w^{\alpha -1}\over w^2 +1}\,dw
\ \ \hbox{for}\ \ \nu\to 0
$$
and with $\beta =    \alpha -1$, hence $-1<\beta <1$, we have to
determine the value of
$$
q(\beta) =  \int_0^\infty {x^\beta\over x^2 +1}\, dx\,.
$$
Observing that $\,q(-\beta)= q(\beta)\,$ (substitute $\xi = 1/x$) we
do this  $\,0\le \beta<1\,.$
Complementation by (integrate
along the upper edge of the negative real semi-axis)
$$
\int_{-\infty}^0 {x^\beta\over x^2 +1}\,dx =  e^{i\beta\pi} \,
\int_0^{+\infty} {x^\beta\over x^2 +1}\,dx
$$
gives, via the residue theorem,
$$
 \l(1+e^{i\beta \pi}\r)\,q(\beta) =	\int_{-\infty}^{+\infty}
{x^\beta\over x^2 +1}\,dx= \pi\, i^\beta  = \pi \,e^{i \beta \pi/2}
\,. $$
So
$$    q(\beta ) = {\pi \over  2\cos{(\beta\pi/ 2)}} =
 {\pi \over  2\sin{(\alpha\pi/ 2)}}  \,,
$$
and hence
$$
\rho(\nu) \rightarrow {\pi \over 2\sin{(\alpha \pi /2)}}
\q \hbox{if}\q 0<\alpha <2 \q (\nu \to 0+).
\leqno(A.4) $$
\vsp
In case (ii) the integral diverges for $\nu =  0$, so we must proceed in
another way. Inserting $\alpha =2$ in (A.3) and
differentiating we obtain for $\nu >0$
$$
 -\rho'(\nu) =\int_0^\infty {w^2 e^{-\nu w}\over w^2+1}\, dw =
 \int_0^\infty e^{-\nu w}\left(1-{1\over w^2 +1}\right) \,dw =
 {1\over \nu} -{\pi\over 2} + o(1)\,,
$$
and then by integration
$$
\rho(\nu) \sim -\log \nu =    \log {1\over \nu}
\q (\nu \to 0+).
\leqno(A.5) $$
\vsp
Now we can collect results. From (A.1) - (A.5),
using $\nu= \kappa h$
which because of symmetry we can replace by $|\kappa|\,h$
(admitting also negative values of $\kappa $) we deduce
$$
\Re \gamma(z) \sim {\pi \over 2\Gamma(\alpha +1)
 \sin{(\alpha \pi/2)}} |\kappa|^\alpha h^\alpha \ \ \hbox{if}\ \
0<\alpha <2,\ \kappa\not =  0\,,  \q \hbox{as} \q h\to 0\,,
\leqno(A.6) $$
$$
\Re \gamma(z) \sim \kappa^2 h^2 \log {1\over |\kappa|h}\ \ \hbox{if}\ \
\alpha =    2,\ \kappa\not =	0\,, \q \hbox{as} \q h\to 0\,.
\leqno(A.7) $$


\subsection*{Acknowledgements}
We are grateful to the Italian
Istituto Nazionale di Alta Matematica and to the
Research Commission of Free University of Berlin for supporting
the joint efforts of our research groups in Berlin and Bologna.

\end{document}